\documentclass[12pt]{article}
\usepackage[english]{babel}
\usepackage{graphics}
\usepackage{amsfonts}
\usepackage{amssymb}
\usepackage{mathrsfs}
\usepackage{enumerate}

\newtheorem{lemma}{Lemma}[section]

\newtheorem{proposition}{Proposition}[section]

\newtheorem{theorem}{Theorem}[section]
\newtheorem{problem}{Problem}[section]
\newtheorem{remark}{Remark}[section]
\newtheorem{definition}{Definition}[section]
\def\QED{\hfill $\Box$\par\smallskip\noindent}

\def\scatola{\lower5pt\hbox{\vbox{\hrule\hbox{\vrule\kern2pt\vbox%
{\kern5pt\hbox{\mathsurround=0pt
}\kern2pt}\kern4pt\vrule}\hrule}}\ } 
\def\l{\lambda}

\def\noin{\noindent}
\def\H{{I\!\!H}}

\def\R{{I\!\!R}}

\def\l{\lambda}

\def\X{\mathbb{X}}
\def\go{\overline g}
\def\x{\overline x}
\def\y{\overline y}
\def\t{\overline t}

\def\0{\overline 0}
\def\pn{\par\noindent}
\def\pmn{\par\medskip\noindent}
\def\psn{\par\smallskip\noindent}

\font\tenmsb=msbm10 \font\sevenmsb=msbm7 \font\fivemsb=msbm5
\newfam\msbfam
\textfont\msbfam=\tenmsb \scriptfont\msbfam=\sevenmsb
\scriptscriptfont\msbfam=\fivemsb
\font\teneufm=eufm10 \font\seveneufm=eufm7 \font\fiveeufm=eufm5
\newfam\eufmfam
\textfont\eufmfam=\teneufm \scriptfont\eufmfam=\seveneufm
\scriptscriptfont\eufmfam=\fiveeufm
\begin{document}
\author{A. Calogero\thanks{Dipartimento di Statistica,
Universit\`a degli Studi di Milano Bicocca, Via Bicocca degli
Arcimboldi 8, I-20126 Milano \tt(andrea.calogero@unimib.it)} and
R. Pini\thanks{Dipartimento di Metodi Quantitativi per le Scienze
Economiche ed Aziendali, Universit\`a degli Studi di Milano
Bicocca, Via Bicocca degli Arcimboldi 8, I-20126 Milano
\tt(rita.pini@unimib.it)}}

\title{Horizontal normal map on the Heisenberg group}

\maketitle
\begin{abstract} \noindent We investigate the notion of H--subdifferential and H--normal
map of a function $u$ on the Heisenberg group, based on the
sub--Riemannian structure of $\H.$ We show that some properties of
the subdifferential in the Euclidean setting are inherited. In
particular, a characterization of the convexity of a function is
given via the nonemptiness of the H--subdifferential
$\partial_Hu(g)$ at every point $g.$ Concerning the H--normal map,
we prove a monotonicity result when suitable strictly convex
radial functions are considered. Finally, we suggest a definition
of the Monge--Amp\`ere measure of a function $u$ via its H--normal
map, and we extend a well--known integration result by
Rockafellar.
\end{abstract}


\section{Introduction}

The notion of subdifferential and normal map for a function defined
on Euclidean spaces, or, more generally, in Banach spaces, is a
classical concept (see, for instance, \cite{Ro1969}, \cite{Lu1969}).
The most interesting features concern the properties of the
subdifferential arising from convex functions. Indeed, in this case,
it enjoys some more interesting properties, among them its
uniqueness at a point characterizes the differentiability of the
function at the same point.

The notion of normal map has been exploited in order to define the
weak solutions to the Monge--Amp\`ere equation $\det D^2u=f$ (see,
for instance, \cite{Gu2001} and the references therein). In
particular, a pointwise estimate for convex functions, the
Alexandrov's maximum principle, is of great importance in the theory
of weak solutions for the Monge--Amp\`ere equation; its proof relies
on a monotonicity property of the normal map of convex functions,
and it is based on geometric features of their graphs.

A natural question rises whether similar properties and maximum
principle results can be stated in the more general setting of
Carnot groups or, more specifically, of the Heisenberg group.

Many papers have been devoted to the study of different types of
convexity in Carnot groups. Seemingly, the most interesting and
fruitful notion of convexity is the so--called weakly H--convexity
(see \cite{DaGaNh2003}, \cite{LuMaStr2004}). In \cite{DaGaNh2003},
in order to deal with the horizontal version of the
Monge--Amp\`ere operator in a Carnot group, the authors worked out
a notion of horizontal subgradient and of horizontal normal map of
a function $u$ that are strictly related to the definition of weak
H--convexity; in the paper, we will refer to weak H--convexity
simply as to {\it convexity.}

In \cite{GuMo2005}, the authors investigated a Monge--Amp\`ere
type operator on the Heisenberg group, but they followed a
different route with respect to the Euclidean framework (see
\cite{Gu2001}), expressing their doubts about the existence of a
suitable definition of normal map in $\H$ useful to state maximum
comparison results.

The aim of this paper is to shed some more light about properties
and potential of the normal map of a function in the Heisenberg
setting.

In Section \ref{basic notions} we provide the main definitions,
while in Section \ref{H-subdifferentiability} we define the
H--subdifferential and the associated H--normal map for a function
$u:\H\to\R,$ that takes into account the sub--Riemannian structure
of the Heisenberg group.

In Section \ref{H--subdifferentiability and convex functions} we
state some results concerning the H--subdifferential of a convex
function on $\H;$ in particular, we show that, for convex functions,
 the H--subdifferential is nonempty at every point.
 Moreover, the uniqueness of the H--subdifferential at a
point $g_0$ is equivalent to the existence of $\X u(g_0).$

The main purpose of Section \ref{The H--normal map} is the
investigation of the normal map $\partial_Hu;$ in particular, we
are interested in those properties of this map that are inherited
from the corresponding properties of the map in Euclidean spaces.
We show that the image of compact subsets of $\H$ under the map
$\partial_Hu$ are compact subsets of $V_1.$ Furthermore, we prove
a monotonicity result for the H--normal map of strictly convex,
radial functions satisfying an additional assumption.

In Section \ref{applications}, in order to show how the
H--subdifferential of a function $u$ carries much information
about the function itself, we suggest a definition of the
Monge--Amp\`ere measure of $u$ via its H--normal map, and we prove
an extension of a well--known integration result due to
Rockafellar.

\section{Basic notions}\label{basic notions}
\vskip0.5truecm

The Heisenberg group $\H=\H^1$ is the Lie group given  by the
underlying manifold $\R^3$ with the non commutative group law
$$
gg'=(x,y,t)(x',y',t')= \left(x+x',y+y',t+t'+2(x'y-xy')\right).
$$
\noin The unit element is $e=(0,0,0),$ and the inverse of
$g=(x,y,t)$ is $g^{-1}=(-x,-y,-t)$. Left translations and
anisotropic dilations are, in this setup, $L_{g_0}(g)=g_0g$ and
$\delta_\lambda(x,y,t)=\left(\lambda x,\lambda y,\lambda^2t\right).$

The  differentiable  structure on $\H$ is determined  by  the left
invariant vector fields
$$
X=\frac{\partial}{\partial x}+2y\frac{\partial}{\partial t},\qquad
Y=\frac{\partial}{\partial y}-2x\frac{\partial}{\partial t},\qquad
T=\frac{\partial}{\partial t},\qquad {\rm with}\ \ [X,Y]=-4T.
$$
The vector field $T$ commutes with the vector fields $X$ and $Y$;
$X$ and $Y$ are called \it horizontal vector fields\rm.

The Lie algebra of $\H$, $\mathfrak{h}$, is the stratified algebra
$\mathfrak{h}=\R^3=V_1\oplus V_2,$ where $V_1={\rm
span}\left\{X,Y\right\},$ $V_2={\rm span}\left\{T\right\};$
$\langle\cdot ,\cdot\rangle$ will denote the inner product. Via
the exponential map exp we identify the vector $\alpha X+\beta
Y+\gamma T$ in $\mathfrak{h}$ with the point $(\alpha, \beta,
\gamma)$ in $\H;$ the inverse $\xi: \H\to\mathfrak{h}$ of the
exponential map has the unique decomposition $\xi=(\xi_1,\xi_2)$
with $\xi_i:\H\to V_i.$

The main issue in the analysis of the Heisenberg group is that the
classical first and second order differential operators are
considered only in terms of horizontal fields. For a given open
subset $\Omega\subset\H$, the class $\Gamma^k(\Omega)$ represents
the Folland--Stein space of functions having continuous
derivatives up to order $k$ with respect to the vector fields $X$
and $Y;$ we denote as usual by ${\cal C}^k(\Omega)$ the class of
functions having continuous derivatives up to order $k$ with
respect to the differential structure of $\R^3$.

Let us recall that the horizontal gradient of a function
$u\in\Gamma^1(\Omega)$ at $g\in\Omega$ is the $2$--vector
$$
(\nabla_h u)(g)= \left((Xu)(g), \ (Yu)(g)\right),
$$
written with respect to the basis $\{X,Y\}$ of $V_1$ ; we denote
by $\X u$ the element in $V_1$ defined as follows
$$
\X u=(Xu)X+(Yu)Y.
$$
The horizontal Hessian of $u\in\Gamma^2(\Omega)$ at $g\in\Omega$ is
the $2\times 2$ matrix
$$
(\nabla_h^2u)(g)= \left(\matrix{(X(Xu))(g)&(X(Yu))(g)\cr
(Y(Xu))(g)&(Y(Yu))(g)\cr} \right),
$$
while the symmetrized horizontal Hessian is the $2\times 2$
symmetric matrix
$$
\left[(\nabla_h^2u)(g)\right]^*= \frac{1}{2}\left\{
(\nabla^2_hu)(g)+\left[(\nabla^2_hu)(g)\right]^T\right\}.
$$

Given a point $g_0\in\H$, the horizontal plane $H_{g_0}$ associated
to $g_0$ is the plane in $\H$ defined by
$$
H_{g_0}=L_{g_0}\left(\exp(V_1)\right)=\left\{g=(x,y,t)\in\H:
t=t_0+2y_0x-2x_0y\right\}.
$$
Notice that, given $g\in\H_{g_0}$, the set
$\{g_0\delta_{\lambda}(g_0^{-1}g), \, \lambda\in [0,1]\}$ is the
segment in $H_{g_0}$ joining $g_0$ with $g$ (i.e., the convex
closure, in the Euclidean sense, of the set $\{g_0, g\}$). We call
it {\it horizontal segment}.

By contrast with Euclidean spaces, where the Euclidean distance is
the most natural choice, in the Heisenberg group several distances
have been introduced for different purposes. However, all of these
distances are homogeneous, namely, they are left invariant and
satisfy the relation $\rho(\delta_r g',\delta_r g)=r\rho(g',g)$ for
every $g',g\in \H,$ and $r>0.$ In particular, the Euclidean distance
to the origin $|\cdot|$ on $\mathfrak{h}$ induces a homogeneous
pseudo--norm $||\cdot||$ on $\mathfrak{h}$ and, via the exponential
map, on the group $\H$ in the following way: for $v\in
\mathfrak{h},$ with $v=v_1+v_2,$ $v_i\in V_i,$ we let
$$
||v||=(|v_1|^4+v_2^2)^{1/4},
$$
and then define the pseudo--norm on $\H$ by the equation
$$
\rho(g)=||v||,\qquad {\mathrm{if}\,}g=\exp v.
$$
The distance between $g$ and $g'$ is given by $\rho(g^{-1}g').$

\vskip0.5truecm
\subsection{Pansu differentiability}
\vskip0.5truecm

Let $\rho$ be any homogeneous distance on $\H$; to simplify the
notation, we denote by $\rho(g)$ the distance $\rho(g,e).$

Let $u:\H\to \R^k.$ We say that $u$ is Pansu differentiable at
$g\in\H$ if there exists a $G$--linear map $Du(g):\H\to\R^k,$ i.e.,
a group homomorphism that satisfies the relation $Du(g)(\delta_r
h)=r Du(g)(h)$ for every $h\in\H$ and $r>0,$ and such that
$$ \lim_{\rho(h)\to
0}\frac{|u(gh)-u(g)-Du(g)(h)|}{\rho(h)}=0.
$$
We call the map $Du(g)$ the Pansu differential of $u$ in $g.$

\noindent The $k\times 3$ matrix representing the Pansu differential
$Du$ of $u=(u^1,u^2,\ldots,u^k)$ can be written as follows
\begin{equation}\label{Pansu2}
    \left(
    \begin{array}{ccc}
    Xu^1&Yu^1&0\\
    Xu^2&Yu^2&0\\
    \ldots&\ldots&\ldots\\
    Xu^k&Yu^k&0
    \end{array}
    \right).
\end{equation}

The horizontal jacobian $J_H u(g)$ of $u$ at $g$ is defined by
taking the standard jacobian of the matrix (\ref{Pansu2}); to this
concern see \cite{EvGa1991}. \noindent In the particular case
$k=1,$ an easy computation gives us that $u$ is Pansu
differentiable at $g$ if
\begin{equation}\label{defPansu}
 Du(g)(h)=\lim_{\lambda\to
0^+}\frac{u(g\delta_\lambda(h))-u(g)}{\lambda}
\end{equation}
exists for every $h\in\H.$ \psn Moreover, if $u\in \Gamma^1(\H),$
the Pansu differential $Du(g)$ is given by the formula
$$
Du(g)(h)=\langle\X u(g),\xi_1(h)  \rangle,
$$
for every $g$ and $h$ in $\H$ (see \cite{DaGaNh2003}).

Let us recall the following relevant definition concerning the
degree of regularity of a function.

\psn \begin{definition} {\rm Let $\Omega\subset \H$ be a bounded
open subset, and $0<\alpha\le 1.$ A bounded function $u:\Omega\to
\R$ is said to belong to the class $\Gamma^{0,\alpha}(\Omega)$ if
there exists a positive constant $L_{\alpha}>0$ such that
$$
|u(g)-u(g')|\le L_{\alpha}\rho(g,g')^{\alpha},\quad g,g'\in \Omega.
$$
A function $u\in \Gamma^{0,1}(\Omega)$ is said to belong to the
class $\Gamma^{1,\alpha}(\Omega)$ if both $Xu$ and $Yu$ exist in
$\Omega$ and $Xu,Yu\in \Gamma^{0,\alpha}(\Omega).$}
\end{definition}

As usual we say that $u$ is Lipschitz continuous if $u\in
\Gamma^{0,1};$ the symbol $\Gamma^{0,1}_{\rm loc}(\Omega)$ denotes
the class of locally Lipschitz continuous functions on $\Omega.$

In the fundamental paper \cite{Pa1989}, Pansu provides a
Rademacher--Stefanov type result in the Carnot group setting; in
particular, he proves that the Lipschitz continuous functions are
differentiable almost everywhere in the horizontal directions. A
further result, due to Danielli, Garofalo and Salsa
(\cite{DaGaSa2003}, Th. 2.7), will play a crucial role for some
results in the sequel. We state it assuming that the Carnot group is
the Heisenberg group.\psn
\begin{theorem}\label{DaGaSa}
Let $\Omega$ be an open subset of $\H,$ and consider $u:\Omega\to
\R,$ with $u\in \Gamma^{0,1}(\Omega).$ Then there exists a set
$E\subset \Omega$ of Haar measure zero such that the Pansu
differential $Du(g)$ and the horizontal gradient $\X u(g)$ exist for
every $g\in \Omega\setminus E,$ and
\begin{equation}\label{Stefanov}
Du(g)(h)=\langle \X u(g),\xi_1(h)\rangle,\qquad \mathrm{for\;
every\; }h\in \H.
\end{equation}
Furthermore, $\X u\in L^{\infty}(\Omega).$
\end{theorem}

\subsection{Convexity}
\vskip0.5truecm

In the Heisenberg group, and in Carnot groups in general, several
definitions of convexity have been introduced and studied for both
sets (see \cite{MoRi2005}, \cite{CaCaPi2007}) and functions (see
\cite{DaGaNh2003}, \cite{LuMaStr2004}). As a matter of fact, the
results obtained in literature suggest that, among them, the most
suitable and satisfactory is the notion of weak H--convexity. In
the sequel, for the sake of brevity, we will refer to weak
H--convexity as to {\it convexity}; to avoid misunderstandings,
the classical convexity will be called {\it Euclidean convexity}.

\begin{definition} {\rm A subset $\Omega$ of $\H$ is said to be {\it convex}
if, for every $g\in \Omega$ and for every $g'\in H_g\cap \Omega,$
$$
g\delta_{\l}(g^{-1}g')\in \Omega,\qquad \forall \l\in [0,1].
$$}
\end{definition}

\begin{definition} {\rm Let $\Omega$ be a convex subset of $\H.$
A function $u:\Omega\to \R$ is called {\it convex} if
 \begin{equation}\label{Def wHconvex}
u(g_\lambda)\le u(g)+\lambda\left(u(g')-u(g)\right)
\end{equation}
for all $g\in\Omega$, $g'\in H_g\cap\Omega,$ and $\lambda\in[0,1].$
This is equivalent to say that
$$
u(g\exp(\l v))\le u(g)+\lambda\left(u(g\exp(v))-u(g)\right)
$$
for every $g\in\Omega,$ $v\in V_1$ and $\l\in [0,1].$

We say that $u$ is strictly convex if is convex and the equality
in (\ref{Def wHconvex}) holds, whenever $g\neq g',$ if and only if
$\l=0$ or $\l=1.$}
\end{definition}

Observe that $u$ is a convex function on $\Omega$ if and only if
$u$ is Euclidean convex on any horizontal segment.

\pn In the sequel, without saying it explicitly, we will assume
that the domain of a convex function is an open convex set.

Next theorem (see, for instance, \cite{DaGaNh2003}) provides a
useful second order condition for convexity, based on the behaviour
of the symmetrized Hessian of $u.$
\begin{theorem}\label{car DGN convex function}
Let $\Omega$ be an open convex subset of $\H$ and let
$u\in\Gamma^2(\Omega).$ Then, $u$ is convex on $\Omega$ if and only
if the symmetrized horizontal Hessian
$\left[(\nabla_h^2u)(g)\right]^*$ is positive semidefinite for every
$g\in\Omega$.
\end{theorem}
In \cite{CaCaPi2008} the authors provide a characterization of
quasi--convex functions in  $\mathcal{C}^2(\H)$ involving the
symmetrized horizontal Hessian as well.

It is worthwhile to mention the following regularity result about
convex functions on $\H$ due to Balogh and Rickly:
\begin{theorem}\label{BaRi} (see \cite{BaRi2003}, Theorem 1.2)
Let $u:\H\to\R$ be a convex function. Then $u\in
\Gamma^{0,1}_{\mathrm{loc}}.$
\end{theorem}

\vskip0.5truecm
\section{H--subdifferential and H--normal map}\label{H-subdifferentiability}
\label{subsec H-sub} \vskip0.5truecm

Let $\Omega'\subset\R^n$ be an open set. Let us recall (see
\cite{Ro1969}) that for every function $f:\Omega'\to\R,$ the {\it
subdifferential of $f$} at a point $x_0$ is defined as follows:
\begin{equation}\label{defEucsubdif}
\partial f(x_0)=\{p\in\R^n:\; f(x)\ge f(x_0)+\langle p,
x-x_0\rangle,\quad \forall x\in\Omega'\}.
\end{equation}
If $\partial f(x_0)$ is not empty, we say that $f$ is {\it
subdifferentiable} at $x_0.$

\noindent The {\it normal map} of $f$ is the set--valued function
$\partial f:\mathcal{P}(\Omega')\to \mathcal{P}(\R^n)$ defined by
\begin{equation}\label{defEucnormap}
\partial f(E)=\bigcup_{x\in E}\partial f(x),
\end{equation}
 for every $E\subset\Omega'.$

In \cite{DaGaNh2003} a notion of horizontal subdifferential that
takes into account the sub--Riemannian structure of $\H$ is given.

\begin{definition}\label{Def
Hsubdiffer1} {\rm Let $u:\Omega\to\R,$ with $\Omega$ open subset
of $\H.$ The {\it horizontal subdifferential} (or {\it
H--subdifferential}) of $u$ at $g_0\in\Omega$ is the set
$$
\partial_Hu(g_0)=\{ p\in V_1:\; u(g)\ge u(g_0)+\langle p,
\xi_1(g)-\xi_1(g_0)\rangle, \quad \forall g\in H_{g_0}\cap\Omega\}.
$$}
\end{definition}
As in the classical context, we say that $u$ is {\it horizontally
subdifferentiable} (shortly, {\it H--subdifferentiable}) at $g_0$
if $\partial u_H(g_0)$ is not empty. Moreover, if $p\in\partial
u_H(g_0),$ we say that $p$ is a {\it H--subgradient} of $u$ at
$g_0.$

In \cite{DaGaNh2003} the authors proved the following result:

\begin{proposition}\label{prop} (see \cite{DaGaNh2003},  Proposition 10.6).
Let $u$ be a function in $\Gamma^1(\Omega),$ and $\Omega\subset \H$
open. If $\partial_Hu(g)\neq \emptyset$, then
$\partial_Hu(g)=\{\mathbb{X}u(g)\}.$
\end{proposition}

\vskip0.5truecm
\subsection{An equivalent notion of H--subdifferentiability}
\vskip0.5truecm

In the Euclidean context another notion of subdifferentiability can
be done and we say (see, for instance, \cite{GaMc1996} and
\cite{GuNg2007}) that $f$ is subdifferentiable at $x_0\in
\Omega'\subset\R^n$ if there exists $p\in\R^n$ such that
\begin{equation}\label{defEucsubdif2}
f(x)\ge f(x_0)+\langle p, x-x_0\rangle +o(|x-x_0|)
 \end{equation}
 as $|x-x_0|\to 0.$
This notion is useful in the study of optimal mass transportation
problems together with the notions of $c$--subdifferentiability,
$c$--convexity 
and Legendre--Fenchel transform.

Recently, these concepts have been investigated in the framework of
the Heisenberg group. For instance, in \cite{AmRi2004} the authors
defined $c$--subdifferentiability and $c$--convexity for functions
on $\H;$ in \cite{CaPi2009}, taking into account the horizontal
structure, the Fenchel transform was introduced for functions on
$\H,$ and a characterization of convexity was provided.

Starting from (\ref{defEucsubdif2}), another notion of
H--subdifferentiability can be given:
\begin{definition}\label{Def Hsubdiffer2}
{\rm  Let $u:\H\to\R.$ We say that $u$ is {\it horizontally
subdifferentiable} at $g_0\in \Omega$ if there exists $p\in V_1$
such that
$$
u(g)\ge u(g_0)+\langle p, \xi_1(g)-\xi_1(g_0)\rangle
+o(||\xi_1(g)-\xi_1(g_0)||), \quad \texttt{\rm for}\ g\in H_{g_0}\
\texttt{\rm and}\ g\to g_0.
$$}
\end{definition}
In next Proposition we show that, in the case of convex functions,
the two notions of H--subdifferentiability given in Definitions
\ref{Def Hsubdiffer1} and \ref{Def Hsubdiffer2} are equivalent.

\begin{proposition}
Let $u$ be a convex function. Then $u$ is H--subdifferentiable at
$g_0$ (in the sense of definition \ref{Def Hsubdiffer2}) if and only
if $\partial u_H(g_0)$ is not empty.
\end{proposition}
{\bf Proof:} If $\partial_Hu(g_0)\neq \emptyset,$ the $u$ is
trivially H--subdifferentiable at $g_0$ according to Definition
\ref{Def Hsubdiffer1}. Assume that there exists $p\in V_1$ such
that definition \ref{Def Hsubdiffer2} is fulfilled at $g_0.$ Let
us prove that $p\in \partial_Hu(g_0).$ By contradiction, let
$g'\in H_{g_0}$ such that
$$
u(g')-u(g_0)-\langle p,\xi_1(g')-\xi_1(g_0)\rangle=\alpha<0.
$$
Define the function $U:[0,1]\to \R$ as follows:
\begin{eqnarray*}\label{a}
U(\l)& = &u(g_{\l})-u(g_0)-\langle
p,\xi_1(g_{\l})-\xi_1(g_0)\rangle\\ & = &u(g_{\l})-u(g_0)-\l\langle
p,\xi_1(g')-\xi_1(g_0)\rangle,
\end{eqnarray*}
where $g_{\l}=g_0\delta_{\l}(g_0^{-1}g')$ varies along the
horizontal segment $[g_0,g']$ as $\l$ varies in $[0,1].$ The
function $U$ is Euclidean convex, $U(0)=0,$ $U(1)=\alpha<0;$ in
particular, for every $\l\in [0,1],$
\begin{eqnarray}\label{U}
U(\l)& \le & (1-\l)U(0)+\l U(1)\\ &=& \l\alpha.
\end{eqnarray}
At the same time, by the assumption of H--differentiability of $u$
at $g_0,$
\begin{equation}\label{b}
U(\l)\ge
o(||\xi_1(g_{\l}-g_0)||)=||\xi_1(g')-\xi_1(g_0)||\,o(|\l|),\qquad
\l\to 0^+.
\end{equation}
Putting together (\ref{U}) and (\ref{b}), we get
$$
||\xi_1(g')-\xi_1(g_0)||o(|\l|)\le \l\alpha,\qquad \l\to 0^+,
$$
or, dividing by $\l,$
$$
||\xi_1(g')-\xi_1(g_0)||o(1)\le \alpha,\qquad \l\to 0^+,
$$
a contradiction, since $\alpha<0.$ \QED \pmn

In the sequel, we prefer to deal with the notion of
H--subdifferentiability of Definition \ref{Def Hsubdiffer1}.

\vskip0.5truecm
\subsection{The Horizontal normal map}
\vskip0.5truecm

The notion of horizontal normal map associated to the horizontal
subdifferential arises naturally:

\begin{definition}{\rm
Let $u:\Omega\to(-\infty,+\infty],$ $\Omega$ open. The {\it
horizontal normal map} of $u$ (or {\it H--normal map}) is the set
valued function $\partial_H u:\mathcal{P}(\Omega)\to
\mathcal{P}(V_1)$ defined by
$$\partial_H u(E)=\bigcup_{g\in E}\partial_H u(g),
$$
for every $E\subset\Omega.$}
\end{definition}

One of the purposes of this paper is to establish whether the
H--normal map can play a suitable role in dealing with the
Monge--Amp\`ere measure of a function (see \cite{GuMo2004}), or
with some maximum or comparison principles for convex functions
(see \cite{GuMo2005}).

Let us recall that, given a map
$\mathcal{F}:\mathcal{P}(\Omega)\to \mathcal{P}(\mathbf{Y}),$
where $\Omega\subset \mathbf{X},$ the graph of $\mathcal{F}$ is
defined as the set
$$
\{(x,y)\in \mathbf{X}\times \mathbf{Y}:\, x\in \Omega,\,y\in
\mathcal{F}(x)\}.
$$
A possible extension of the concept of continuity of a function to
maps is provided by the notion of closed graph.

\begin{definition}{\rm
Let $\mathbf{X},\mathbf{Y}$ be topological spaces. A map
$\mathcal{F}:\mathbf{X}\to \mathcal{P}({\mathbf{Y}})$ is said to
have closed graph if for every $\{x_n\}\subset \mathbf{X},$
$x_n\to x\in \mathbf{X},$ $\{y_n\}\subset \mathbf{Y},$ with
$y_n\in \textbf{F}(x_n),$ then
$$
y_n\to y\Longrightarrow y\in \textbf{F}(x).
$$}
\end{definition}

\vskip0.5truecm
\section{H--subdifferentiability and convex functions}\label{H--subdifferentiability and convex functions}
\vskip0.5truecm

The aim of this section is to investigate some properties of the
H--subdifferential and of the H--normal map of convex functions on
$\H.$ Our main result is Theorem \ref{teorema mappa normale},
where a convex function is characterized via its
H--subdifferentiability on the domain. As a consequence, the
uniqueness of the H--subdifferential of a convex function $u$ at a
point $g$ is equivalent to the existence of $\X u(g).$

If $u:\Omega\to\R$ is convex, then, for every $g\in \Omega,$ the
limit
$$
\lim_{\l\to 0^+}\frac{u(g\exp(\l v))-u(g)}{\l}
$$
exists in $\mathbb{R},$ for every $v\in V_1.$ We set
$$
u'(g;v)=\lim_{\l\to 0^+}\frac{u(g\exp(\l v))-u(g)}{\l}.
$$

Let us first state the following useful characterization of
H--subdifferentials of a convex function.
\begin{proposition}\label{caratterizzazione sottogradiente}
Let $u:\Omega\to \R$ be a convex function, and $g\in \Omega.$ Then
$p\in\partial_Hu(g)$ if and only if
\begin{equation}\label{prop diff enunc}
u'(g;v)\ge \langle p,v\rangle,\qquad \texttt{\rm for every}\ v\in
V_1.
\end{equation}
In particular, if $u$ is Pansu differentiable at $g$ and
(\ref{Stefanov}) holds, then
$$
\X u(g)\in \partial_Hu(g).
$$
\end{proposition}
{\bf Proof:} Suppose that $p\in\partial_Hu(g).$ Hence, for every
$\l\in [0,1]$ and $v\in V_1,$ we have
 \begin{eqnarray*}
 u(g\exp (\l v))
 &\ge& u(g)+\l\langle p,v\rangle\label{prop diff}.
 \end{eqnarray*}
The previous inequality (\ref{prop diff}) holds if and only if
$$
\frac{ u(g\exp (\l v))- u(g)}{\l} \ge \langle p,v\rangle,
$$
for every $\l\in (0,1]$ and $v\in V_1:$ hence (\ref{prop diff
enunc}) follows obviously.

Conversely, assume that (\ref{prop diff enunc}) holds. Since $u$ is
convex, we have
\begin{equation}\label{prop diff 1}
u(g\exp(\l v))\le u(g)+\l(u(g\exp (v))-u(g)),
\end{equation}
for  every $\l\in [0,1]$ and $v\in V_1.$ By (\ref{prop diff enunc})
and (\ref{prop diff 1}) we get that
$$
\langle p,v\rangle\le \lim_{\l\to 0^+}\frac{u(g\exp(\l
v))-u(g)}{\l}\le u(g\exp (v))-u(g),
$$
for every $v\in V_1;$ hence $p\in\partial_Hu(g).$

Let us suppose now that $u$ is Pansu differentiable at $g$ and
(\ref{Stefanov}) holds; then, by (\ref{defPansu}),
$$
\langle \X u(g),\xi_1(h)\rangle=Du(g)(h)=\lim_{\l\to
0^+}\frac{u(g\delta_\l(h))-u(g)}{\l},
$$
for every $h\in\H.$ The convexity of $u$ gives us that
$\frac{u(g\delta_\l(h))-u(g)}{\l}$ decreases when $\l\to 0^+,$ for
every fixed $g$ and $h\in H_e.$ Hence
$$
\langle \X u(g),\xi_1(h)\rangle\le
\frac{u(g\delta_\l(h))-u(g)}{\l},\quad \forall \lambda\in (0,1],
$$
so that $\X u(g)$ is a H--subdifferential. \QED \pmn

In \cite{DaGaNh2003}, Danielli, Garofalo and Nhieu proved the
following
\begin{proposition}
\label{prop0} (see \cite{DaGaNh2003},  Proposition 10.5). Let
$u:\Omega\to\R,$ where $\Omega$ is an open and convex subset of
$\H.$ If $\partial_Hu(g)\neq \emptyset$ for every $g\in \Omega,$
then $u$ is convex.
\end{proposition}
In order to show that the converse holds, we need next result.
Indeed, the following Lemma will be crucial also in the proof of
Theorem \ref{compatti in compatti}.

\begin{lemma}\label{lemma mappa normale}
Let $\Omega$ be an open subset of $\H,$ and consider a function
$u\in \mathcal{C}(\Omega).$ Then the map
$\partial_Hu:\mathcal{P}(\Omega) \to \mathcal{P}(V_1)$ has closed
graph.
\end{lemma}
\psn {\bf Proof}: We prove the lemma by contradiction. Assume that
$g_0\in \Omega,$ and there exist sequences $\{g_n\}\subset \Omega$
and $\{p_n\},$ with $p_n\in
\partial_Hu(g_n),$ such that
$$
g_n\to g_0,\quad p_n\to p, \quad p\notin
\partial_Hu(g_0).
$$
Consequently, there exists $g'\in H_{g_0}\cap \Omega$ such that
$$
u(g')-u(g_0)=\langle p,\xi_1(g')-\xi_1(g_0)\rangle -\alpha,
$$
for a suitable $\alpha >0.$ From the assumptions,
$u$ is continuous.

Denote by $B(e,r)$ the set of points $g\in \H$ such that
$\rho(g)<r.$ Let us consider the balls $B(e,r)$ and $B(e,r')$ in
such a way that
\begin{enumerate}[i)]\label{bolle}
\item $|u(g'h')-u(g')|<\alpha/10,$ for every $h'\in B(e,r');$
\item $|u(g_0h)-u(g_0)|<\alpha/10,$ for every $h\in B(e,r);$
\item $H_{g_0h}\cap L_{g'}(B(e,r'))\neq \emptyset,$ for every $h\in
B(e,r).$
\end{enumerate}
The reader can easily convince himself that such balls exist using
suitable continuity arguments on $u$ and on the displacement of
the horizontal planes of points moving close to others. For every
$h\in B(e,r)$ and $h'\in B(e,r'),$ i) and ii) imply the following
inequality:
\begin{eqnarray}\label{inequality}
\langle p,\xi_1(g')-\xi_1(g_0)\rangle -\alpha&=&u(g')-u(g_0)\nonumber\\
&\ge & u(g'h')-u(g_0h)-\alpha/5.
\end{eqnarray}
Take $N$ such that $g_n\in L_{g_0}(B(e,r))$ for every $n\ge N,$ and
denote by $h_n$ the point in $B(e,r)$ such that $g_n=g_0h_n.$ Notice
that $h_n\to e.$ Moreover, from the choice of the balls, there
exists $h'_n\in B(e,r')$ such that
$$
g'_n=g'h'_n\in H_{g_0 h_n},\quad g'_n\to g'.
$$
Then, taking into account the assumptions and (\ref{inequality}), we
get
\begin{eqnarray*}
\langle p,\xi_1(g')-\xi_1(g_0)\rangle -\alpha
&\ge &u(g'h'_n)-u(g_0h_n)-\alpha/5\\
&\ge & \langle p_n, \xi_1(g'h'_n)-\xi_1(g_0h_n)\rangle -\alpha/5\\
&\ge & \liminf \langle p_n, \xi_1(g'h'_n)-\xi_1(g_0h_n)\rangle -\alpha/5\\
&=& \langle p,\xi_1(g')-\xi_1(g_0)\rangle -\alpha/5,
\end{eqnarray*}
therefore obtaining
$$
\langle p,\xi_1(g')-\xi_1(g_0)\rangle -\alpha \ge  \langle
p,\xi_1(g')-\xi_1(g_0)\rangle -\alpha/5,
$$
a contradiction. \QED \pmn

We are now able to prove the following interesting
characterization:

\begin{theorem}\label{teorema mappa normale}
Let $u:\Omega\subset \H\to\R,$ where $\Omega$ is open and convex.
Then $u$ is convex  if and only if, for every $g\in \H,$
$\partial_Hu(g)\neq \emptyset.$
\end{theorem}
{\bf Proof}: The \lq\lq if" part is the result of Proposition
\ref{prop0}.

From Theorem \ref{BaRi} any convex function $u$ belongs to
$\Gamma^{0,1}_{loc}(\Omega)$; in particular, it is continuous. By
contradiction, we assume that there exists $g_0\in \Omega$ such
that $\partial_Hu(g_0)=\emptyset.$ Let us consider a neighborhood
$B(e,r)$ of the origin such that $u$ belongs to
$\Gamma^{0,1}(L_{g_0}(B(e,r))).$ From Theorem \ref{DaGaSa}, there
exists a subset $E$ of $L_{g_0}(B(e,r))$ with null measure such
that for every $g\in L_{g_0}(B(e,r))\setminus E$ there exists the
Pansu differential $Du(g)$ and (\ref{Stefanov}) holds; Proposition
\ref{caratterizzazione sottogradiente} shows that, for such $g,$
$\X u(g)\in
\partial_Hu(g).$
Since $\nabla_Hu\in L^\infty(L_{g_0}(B(e,r)))$ (see
\cite{DaGaNh2003}, Theorem 9.1),  there exists $k$ such that $||\X
u(g)||\le k,$ for a.e. $g\in L_{g_0}(B(e,r)).$ Therefore there
exists a sequence $\{g_n\}$ in $L_{g_0}(B(e,r))\setminus E$ such
that
$$
g_n\to g_0,\qquad \X u(g_n)\to p,\qquad \X u(g_n)\in
\partial_Hu(g_n),
$$
for some $p\in V_1.$ Then, since $\partial_H u$ has closed graph
(see Lemma \ref{lemma mappa normale}), $p\in
\partial_Hu(g_0),$ a contradiction.
\QED \pmn

As a matter of fact, the uniqueness of the H--subdifferential for
a convex function characterizes the H--differentiability. We are
able to state the following
\begin{theorem}\label{teorema unicità sottogradiente}
Let $u$ be a convex function on $\Omega.$ Then, $\X u(g)$ exists for
some $g\in \Omega$ if and only if $u$ has a unique H--subgradient at
$g.$ Moreover, in both cases, we have that $\partial_H u(g)=\{\X
u(g)\}.$
\end{theorem}
{\bf Proof:} Let $u$ be convex. By Theorem \ref{teorema mappa
normale}, there exists $p=p_1X+p_2Y\in
\partial_H u(g).$ By Proposition \ref{caratterizzazione
sottogradiente}, if we put $h=(1,0,0)\in H_e,$ we get
\begin{eqnarray*}
p_1&=&\langle p,\xi_1(h)\rangle\\
&\le & \lim_{\l\to 0^+}\frac{u(g\delta_\l(h))-u(g)}{\l}\\
 &=&\lim_{\l\to 0^+}\frac{u(g\exp(\l X))-u(g)}{\l}.
\end{eqnarray*}
Now, taking $h=(-1,0,0)\in H_e,$ similar computations give $$
-p_1\le -\lim_{\l\to 0^-}\frac{u(g\exp(\l X))-u(g)}{\l}.
$$
Hence, we have
$$
p_1=\lim_{\l\to 0}\frac{u(g\exp(\l X))-u(g)}{\l}=Xu(g).
$$
Similar arguments show that $p_2=Yu(g).$

Conversely, assume that $\partial_Hu(g_0)=\{p\}.$ Let us suppose
$\Omega=\H,$ for the sake of simplicity. Fix $v\in V_1,$ and
consider the linear space $\mathcal{L}(v)=\{av,\,a\in\R\}.$ Define
on $\mathcal{L}(v)$ the linear functional
$$
L_v(w)=au'(g_0;v),
$$
for $w=av.$ Notice that $L_v(w)=u'(g_0;w)$ whenever $w=av,$ with
$a>0.$ Indeed,
$$
u'(g_0;w)=\lim_{\l\to 0^+}\frac{u(g_0\exp(\l
w))-u(g_0)}{\l}=\lim_{\l\to 0^+}a\frac{u(g_0\exp(\l
v))-u(g_0)}{\l}=L_v(w).
$$
By the convexity of $u,$ the function $t\to u(g_0\exp(tv))$ is
Euclidean convex; in particular, the following inequality holds
$$
u'(g_0;-v)\ge -u'(g_0;v).
$$
Assume that $w=av,$ with $a<0.$ Then
$$
u'(g_0;w)=u'(g_0; av)=(-a)u'(g_0;-v)\ge au'(g_0;v)=L_v(w).
$$
Since the linear functional $L_v$ satisfies on $\mathcal{L}(v)$ the
inequality
$$
L_v(w)\le u'(g_0;w),
$$
by the Hahn--Banach theorem there exists $p_v\in V_1$ such that
$$
\langle p_v,w\rangle \le u'(g_0;w),\qquad \forall w\in V_1.
$$
From Proposition \ref{caratterizzazione sottogradiente}, $p_v\in
\partial_Hu(g_0),$ and, by the assumptions, $p_v=p.$ In particular, $p_v$ is independent on $v.$
Since $\langle p,v\rangle =u'(g_0;v),$ and $v$ is any vector in
$V_1,$ we can conclude that
$$
\langle p,v\rangle =u'(g_0;v),\qquad \forall v\in V_1,
$$
thus $u'(g_0;\cdot)$ is linear on $V_1.$ Hence, if $w=X$ then
$$
u'(g_0;w)=-u'(g_0;-w)=\lim_{\l\to 0^-}\frac{u(g_0\exp(\l
X))-u(g_0)}{\l}.
$$
This implies the existence of $Xu(g_0),$ that should be equal to
$u'(g_0;X).$ Similar arguments prove the existence of $Yu(g_0);$
in particular, $p=Xu(g_0)X+Yu(g_0)Y.$ \QED

\vskip0.5truecm
\section{The H--normal map}\label{The H--normal map}
\vskip0.5truecm

The aim of this section is to investigate some properties of those
subsets of $V_1$ that are images, via the H--normal map, of subsets
of $\Omega.$ We start by studying the properties of the H--normal
map of a single point, i.e. the H--subdifferential.
\begin{proposition}\label{prop sub compatto}
Let $u:\Omega\to\R,$ with $\Omega\subset\H$ open; let $g\in\Omega.$
Then $\partial u_H(g)$ is a convex set. Moreover, if $u$ is locally
bounded, then $\partial u_H(g)$ is compact.
\end{proposition}
{\bf Proof:} Take any $p_1,p_2\in\partial_Hu(g).$ For every
$\lambda\in [0,1]$ we have that
\begin{eqnarray*}
u(g')&=&(1-\l)u(g')+\l u(g')\\
&\ge& (1-\l)\Bigl(u(g)+\langle p_1, \xi_1(g')-\xi_1(g)\rangle\Bigr)+
\l\Bigl(u(g)+\langle p_2, \xi_1(g')-\xi_1(g)\rangle\Bigr)\\
&=& u(g)+\langle (1-\l)p_1+\l p_2, \xi_1(g')-\xi_1(g)\rangle,
\end{eqnarray*}
for every $g'\in H_g,$ therefore showing that $(1-\l)p_1+\l p_2$
is a H--subgradient of $u.$

Let $\{p_k\}$ be a sequence in $\partial u_H(g).$ For every $k$
and for every $w\in V_1,$ with $||w||=1,$ we have
$$
u(g\delta_\l(\exp(w)))\ge u(g)+\langle p_k,
\xi_1(g\delta_\l(\exp(w)))-\xi_1(g)\rangle= u(g)+\l\langle p_k,
w\rangle.
$$
In particular for every $k$ with $p_k\not=0,$ if we put
$w={p_k}/{||p_k||},$ then we obtain
\begin{equation}\label{compact}
\sup_{||w||=1} u(g\delta_\l(\exp(w)))\ge
u(g)+\frac{\l}{||p_k||}\langle p_k, p_k\rangle= u(g)+\l ||p_k||.
\end{equation}
Take any $\l$ sufficiently small such that
$g\delta_\l(\exp(w))\in\Omega$ for every $w$ with $||w||=1:$ since
$u$ is locally bounded, (\ref{compact}) gives us that $\{p_k\}$ in a
bounded subset of $V_1$. Hence there exists a convergent subsequence
$\{p_{k_n}\}$ such that $p_{k_n}\to p.$ Clearly
$$
u(g')\ge u(g)+\langle p_{k_n}, \xi_1(g')-\xi_1(g)\rangle,
$$
for every $n$ and for every $g'\in\H_g\cap\Omega.$ Letting
$n\to\infty,$ we obtain that $p\in\partial_H u(g).$ \QED \pmn

The investigation of the images of the H--normal map is actually
more awkward. Indeed, if we shift from a point $g$ to another
point $g',$ and consider $p\in\partial_H u(g)$ and
$p'\in\partial_H u(g'),$ then the H--subdifferentials $p$ and $p'$
support the function on the different planes $H_g$ and $H_{g'}.$

First of all, it easy to see that, as in the Euclidean case (see
\cite{Gu2001}), if $K$ is Euclidean convex, then $\partial_H u(K)$
is not necessary Euclidean convex in $V_1.$

To obtain information about $\partial u_H(\Omega),$ we need some
regularity assumptions on $u.$

\begin{proposition}\label{prop sub limitato}
Let $u:\Omega\to\R$ and $u\in\Gamma^{0,1}(\Omega),$ with Lipschitz
constant $L.$ Then $||p||\le L,$ for every $p\in \partial
u_H(\Omega);$ in particular, $\partial u_H(\Omega)$ is a bounded
set.
\end{proposition}
{\bf Proof:} Let $p\in\partial u_H(\Omega).$ Then there exists
$g\in\Omega$ such that, for every $h\in H_e$ e $\l>0,$ we have
$$\frac{u(g\delta_\l(h))-u(g)}{\l}\ge \langle p, \xi_1(h)\rangle.$$
Since $u$ is Lipschitz continuous, there exists $L>0$ such that
$|u(g)-u(g')|\le L \rho(g,g'),$ for every $g'\in\Omega.$ Hence
$$\frac{u(g\delta_\l(h))-u(g)}{\l}\le L \rho(h).$$
The previous two inequalities give us that
$$
\langle p, \xi_1(h)\rangle\le L \rho(h),$$ for every $h\in H_e.$
If we put $h=\xi_1^{-1}(p),$ we obtain $||p||^2\le L ||p||.$
 \QED
\pmn

\begin{theorem}\label{compatti in compatti}
Let $u:\Omega\to\R$ be a function in $\Gamma^{0,1}_{{\rm
loc}}(\Omega).$ Then, for every compact set $K\subset\Omega,$ the
set $\partial u_H(K)$ is compact.
\end{theorem}
{\bf Proof:} By the assumption, for every $g\in\Omega$ there
exists a neighborhood $B_{g}$ such that $u\in \Gamma^{0,1}(B_g),$
i.e., there is a constant $L_{g}$ such that
$$
|u(g'')-u(g')|\le L_{g} \rho(g'',g'),
$$
for every $g'',g'\in B_{g}.$

Let $K$ be any compact subset of $\Omega;$ then
$K\subset\cup_{i=1}^N B_{g_i}$ for a suitable finite set of points
$\{g_1,g_2,\dots,g_N\}.$

Take any $p\in \partial u_H(K),$ and denote by $g$ a point in $K$
such that $p\in \partial u_H(g)\subset \partial_Hu(B_{g_i}),$ for
some $i.$ Hence, by Proposition \ref{prop sub limitato}, since
$u\in\Gamma^{0,1}(B_{g_i})$ for every $i=1,2,\dots, N,$ we have
that
$$
||p||\le L_{g_i}\le \max\{L_{g_i}:\ 1\le i\le N\}.$$
Hence
$\partial u_H(K)$ is bounded.

Let us now consider a sequence $\{p_n\}\subset \partial u_H(K),$
and assume that $p_n\to p,$ and denote by $g_n$ a point in $K$
such that $p_n\in\partial_Hu(g_n).$ Since $K$ in compact, there
exists a subsequence $\{g_{n_k}\}$ such that $g_{n_k}\to g$ for
some $g\in K.$ Since the normal map $\partial_H u$ has closed
graph (see lemma \ref{lemma mappa normale}), $p\in
\partial_H u(g).$
 \QED
Next result follows trivially from Theorems \ref{compatti in
compatti} and \ref{BaRi}.
\begin{remark} {\rm Let $u$ be a convex function on
$\Omega.$ Then, for every compact set $K\subset\Omega,$ the set
$\partial u_H(K)$ is compact.}
\end{remark}

\vskip0.5truecm
\subsection{Monotonicity property of the H--normal map}
\vskip0.5truecm

The purpose of this subsection is to investigate whether, as in the
Euclidean context, a monotonicity property of this type holds:

\medskip
\begin{problem}\label{problem} {\rm Let
$\Omega\subset \H$ be open and bounded, and denote by $u$ and $v$
two convex functions such that
$$
u(g)=v(g),\qquad\forall g\in\partial \Omega $$ and
$$
u(g)\le v(g),\qquad \forall g\in \Omega.
 $$
Then
$$\partial_H v(\Omega)\subset  \partial_H u(\Omega).$$}
\end{problem}

\noindent We are not able to prove the result in the general case;
at the moment, several difficulties rise in the proof. In the sequel
we will consider a particular case of the situation described in
Problem (\ref{problem}); we assume to work with functions
$u:\H\to\R$ with the following property:
$$
u(x,y,t)=U(r,t),
$$
for every $(x,y,t)$ where $r=\sqrt{x^2+y^2}.$ By abuse of language,
we will call these functions \lq\lq radial".

\pn We think that the ideas contained in the proof of this special
case could be of some interest to deal with the general one.

For every $t\in\R$ and $R\ge 0,$  we will denote by $C(t,R)$ the
set
$$
C(t,R)=\{(x,y,t)\in\H:\,x^2+y^2=R^2\},
$$
and by $D(t,R)$ the \lq\lq open\rq\rq\ disc in $H_{(0,0,t)}$
defined by
$$
D(t,R)=\{(x,y,t)\in\H:\,x^2+y^2<R^2\}.
$$
In order to prove our monotonicity result, we need to state the
following technical propositions, that explains some geometric
features of the normal map of radial functions.

\begin{proposition}\label{image of radial is radial}
Let $u:\H\to\R$ be a radial, convex function in $\mathcal{C}^1(\H).$
Then, for every $t\in\R$ and $R\ge 0$,
$$
\partial_Hu(C(t,R))=\{p\in V_1:\, ||p||=R'\},
$$
for a suitable $R'\ge 0.$
\end{proposition}
{\bf Proof:} Under the assumptions, for every $g\in \H$ the map
$\partial_Hu:\H\to V_1$ is single--valued, and $\partial_Hu(g)=\{\X
u(g)\}$ (see Theorem \ref{teorema unicità sottogradiente}).

\pn For every $g=(r\cos\theta,r\sin \theta,t)\in\H,$ we have
$$
Xu(g)=U_r(r,t)\cos\theta+rU_t(r,t)\sin\theta,\quad
Yu(g)=U_r(r,t)\sin\theta-rU_t(r,t)\cos\theta,
$$
and
$$
||\X u(g)||=\sqrt{U_r^2(r,t)+4r^2U_t^2(r,t)}.
$$
If $\theta$ varies in $[0,2\pi)$ we get the thesis. \QED \pmn

\begin{proposition}\label{image of omega} Let $u:\H\to\R$ be a radial, convex
function in $\mathcal{C}^1(\H).$ Denote by $\Omega$ the (nonempty)
sublevel set
$$
\Omega=\{g\in\H:\, u(g)<0\}.
$$
Then, $\partial_Hu(\Omega\cap H_{(0,0,t)})$ is a disc (open or
closed) in $V_1,$ centered at the origin. Consequently, the set
$\partial_Hu(\Omega)$ is a disc (open or closed) in $V_1,$ centered
at the origin.
\end{proposition}
{\bf Proof}: Notice that
$$
\partial_Hu(\Omega)=\bigcup_{t\in\R}\partial_Hu(\Omega\cap
H_{(0,0,t)}),
$$
and $\Omega\cap H_{(0,0,t)}=D(t,R(t)).$ By Proposition \ref{image
of radial is radial} and by the continuity of $\X u$ we easily get
the thesis. \QED \pmn

We would like to emphasize that, in general, $||\X
u(r\cos\theta,r\sin\theta,t)||$ is not an increasing function of
$r;$ this explains why, without further conditions on $u,$ one
cannot infer that $\partial_H u(\partial D(t,R))=
\partial(\partial_H u(D(t,R))).$
The radius $\tilde R$ of $\partial_H u(D(t,R))$ is given by the
expression
$$
\tilde R(t,R)=\sup\{||\X u(r\cos\theta,r\sin\theta,t)||,\,0\le r<
R,\, \theta\in [0,2\pi)\}.
$$

\begin{theorem}
Let $u\in \mathcal{C}^1(\H),$ $v\in \mathcal{C}^2(\H)$ be radial,
strictly convex functions such that $u\le v.$ Denote by $\Omega$
the (nonempty) level set
$$
\Omega=\{g\in \H:\, u(g)<0\}=\{g\in \H:\, v(g)<0\};
$$
assume that $\Omega$ is bounded and
$$
\partial \Omega=\{g\in\H:\, u(g)=0\}=\{g\in\H:\, v(g)=0\}.
$$
If $\go\in\partial \Omega,$ then there exists
$\overline{s}=\overline{s}(\go)\in (0,1]$ such that
$$
\X v(\go)=\overline{s}\X u(\go).
$$
In addition, suppose that the function $V$ defined as
$V(r,t)=v(x,y,t)$ satisfies in $\overline{\Omega}$ the assumption
\begin{equation}\label{condizione monotonia norma}
r^3V_{tr}^2-V_rV_{rr}< 0.
\end{equation}
Then,
\begin{equation}\label{inclusione}
\partial_Hv(\Omega)\subset
\partial_Hu(\Omega).
\end{equation}
\end{theorem}
{\bf Proof:} Take any point $\overline{g}=(\x,\y,\t),$ such that
$\overline{g}\in
\partial\Omega\cap H_{(0,0,\t)}.$ Notice that $\X u(\go)\neq \mathbf{0},$ since $u$ is strictly convex and $\X u(0,0,\t)=\mathbf{0}.$

Let us consider the function $F^u_{\go}:\R^2\to \R$ obtained by
restricting $u$ to $H_{\go}$ and defined as follows:
$$
F^u_{\go}(x,y)=u(x,y,\t+2\y x-2\x y).
$$
Denote by $\Omega^u_{\go}$ the sublevel set
$$
\Omega^u_{\go}=\{(x,y)\in\R^2:\, F^u_{\go}(x,y)< 0\},
$$
and consider, in particular, its boundary
$$
\partial \Omega^u_{\go}=\{(x,y)\in\R^2:\, F^u_{\go}(x,y)= 0\}.
$$
This set is not empty, since $F^u_{\go}(\x,\y)=0;$ moreover, from
the inequality
$$
\X u(\go)=\frac{\partial F^u_{\go}}{\partial
x}(\x,\y)X+\frac{\partial F_{\go}}{\partial y}(\x,\y)Y\neq
\mathbf{0},
$$
the implicit function theorem assures that, at least locally, there
exists a unique curve $\gamma_{\go}:I\to \R^2,$
$\gamma_{\go}(s)=(x(s),y(s)),$ with $(x(s),y(s))\in
\mathcal{C}^1(I)$ and $F^u_{\go}(x(s),y(s))=0;$ moreover, $\nabla
F^u_{\go}(\x,\y)$ is orthogonal to $\dot\gamma_{\go}(0).$ Since
$\nabla F^u_{\go}(\x,\y)$ represents the increasing direction of
$F^u_{\go}$ at the point $(\x,\y),$ the vector $-\nabla
F^u_{\go}(\x,\y)$ points towards $\Omega^u_{\go};$ this implies that
$((\x,\y)-z\nabla F^u_{\go}(\x,\y))$ belongs to $\Omega^u_{\go},$ at
least for small values of $z>0.$

Let us start from the function $v$ instead of $u.$ By the same
arguments applied to $F^v_{\go}=v(x,y,\t+2\y x-2\x y),$ taking into
account that, from the assumptions on $u$ and $v,$ we have
$\Omega^u_{\go}=\Omega^v_{\go},$ we find out that $\nabla
F^u_{\go}(\x,\y)$ and $\nabla F^v_{\go}(\x,\y)$ should satisfy the
equality

$$
\nabla F^v_{\go}(\x,\y)=\overline{s}\nabla F^u_{\go}(\x,\y),
$$
for some positive $\overline{s};$ in other words,

\begin{equation}\label{paralleli su bordi}
\X v(\go)=\overline{s}\X u(\go).
\end{equation}
Let us prove that $\overline{s}\le 1.$ Consider the functions
$f^u,f^v$ defined in $[0,\epsilon)$ for a suitable small $\epsilon$
as follows
$$
f^u(z)=F^u_{\go}((\x,\y)-z\nabla F^u_{\go}(\x,\y)),\quad
f^v(z)=F^v_{\go}((\x,\y)-z\nabla F^v_{\go}(\x,\y)).
$$
Since $f^u(0)=f^v(0)=0,$ and $f^u(z)\le f^v(z)\le 0$ if
$0<z<\epsilon,$ standard arguments of real analysis entail that
$(f^u)'(0)\le (f^v)'(0)\le 0.$ From $(f^u)'(0)=-||\X u(\go)||^2$ and
$(f^v)'(0)=-||\X v(\go)||^2=-\overline{s}^2||\X u(\go)||^2,$  we
obtain that $\overline{s}\le 1.$

Let us consider the function
$$r\mapsto ||\X
v(g)||^2=V_r^2(r,t)+4r^2V_t^2(r,t),
$$
where $g=(r\cos \theta,r\sin\theta, t).$ By the assumption
(\ref{condizione monotonia norma}), standard computations imply that
this is an increasing function; in particular, from Propositions
\ref{image of radial is radial} and \ref{image of omega},
$$
\partial_Hv(\partial (\Omega \cap
H_{(0,0,t)}))=\partial(\partial_Hv(\Omega\cap H_{(0,0,t)})).
$$
Moreover, taking into account that (\ref{paralleli su bordi}) holds
with $\overline{s}\le 1,$ we get
$$
\partial_H v(\Omega\cap H_{(0,0,t)})\subset
\partial_H u(\Omega\cap H_{(0,0,t)}).
$$
Suppose now that $p$ is in $\partial_Hv(\Omega),$ i.e., $p=\X v(g')$
for some $g'=(r'\cos\theta,r'\sin\theta, t')\in \Omega.$ Then
$$
p\in\partial_Hv(C(t',r'))\subset \partial_Hv(\overline{\Omega}\cap
H_{(0,0,t')}) \subset \partial_Hu(\overline{\Omega}\cap
H_{(0,0,t')}),
$$
thereby proving (\ref{inclusione}). \QED \pmn

Following the idea in \cite{BaRi2003}, we consider radial
functions of the type
$$
v(x,y,t)=((x^2+y^2)^2+z(t))^{1/4},
$$
where $z:\R\to\R$ is assumed to be twice continuously
differentiable and positive. Theorem \ref{car DGN convex function}
and easy computations (see \cite{CaCaPi2008}) show that $u$ is
convex on $\H$ if and only if
$$
4z(1+z'')\ge 3(z')^2, \quad \texttt{\rm on}\; \R.
$$
Condition (\ref{condizione monotonia norma}) is equivalent to the
inequality
$$
16(z(t))^2+r^4(16z(t)-(z'(t))^2)\ge 0.
$$

\vskip0.5truecm
\section{Applications}\label{applications}
\vskip0.5truecm

The aim of this section is to show that, like in the Euclidean
framework, the H--subdifferential of a function on the Heisenberg
group carries a lot of information about the function itself.

\vskip0.5truecm
\subsection{The Monge--Amp\`ere measure and H--normal map}
\vskip0.5truecm

In the Heisenberg group, the Monge-Amp\`ere type operator ${\cal
S}_{ma}$ (see \cite{GaTo2005} and  \cite{GuMo2005}) is a fully
nonlinear operator on $u$ defined by
$${\cal S}_{ma}(u)=\texttt{\rm det} [\nabla_H^2 u]^* + \frac{3}{4} ([X,Y]u)^2=
\texttt{\rm det} [\nabla_H^2 u]^* + 12 (Tu)^2
$$
In \cite{GuMo2005} the authors proved the following
result:

\begin{theorem} Given a convex function $u\in \mathcal{C}(\Omega),$
there exists a unique Borel measure $\mu_u$ such that, when $u\in
C^2(\Omega),$
$$
\mu_u(E)=\int_E [{\cal S}_{ma}(u)](g) d g,
$$
for any Borel set $E\subset \Omega.$
\end{theorem}
We call $\mu_u$ the {Monge--Amp\`ere measure} of $u.$

\medskip
In the Euclidean context (see \cite{Gu2001}), the Monge--Amp\`ere
measure $M_f$ associated to a function $f$ is defined via the notion
of normal map $\partial f$ of $f$ (see (\ref{defEucsubdif}) and
(\ref{defEucnormap})). In particular, if $f\in{\cal C}(\Omega'),$
with $\Omega'$ open in $\R^n,$ the Monge--Amp\`ere measure is the
set function $M_f:{\cal E}'\to [0,\infty]$ defined by
\begin{equation}\label{Mf}
M_f(E')=|\partial f(E')|=\int_{\partial f(E')}1\, dp,\qquad \forall
E'\in {\cal E}',
 \end{equation}
where ${\cal E}'=\{E'\subset\Omega':\ \partial f(E')\ \texttt{\rm is
Lebesgue measurable}\}$ and $|A|$ denotes the Lebesgue measure of
$A.$ If $f$ is an Euclidean convex function in $C^2(\Omega'),$ we
have that
\begin{equation}\label{misura MonAmp convex classica}
M_f(E')=\int_{E'} \texttt{\rm det}  [D^2 f](x) d x,
\end{equation}
for every Borel set $E'\subset\Omega'.$ The proof of (\ref{misura
MonAmp convex classica}) (see \cite{Gu2001}) exploits the property
that if $f$ is Euclidean convex and ${\cal C}^2(\Omega'),$ we can
identify $\partial f$ with $\nabla f$ and $\nabla f$ is
one--to--one on the set $\{x\in\Omega':\ D^2f(x)>0\}.$ Hence every
point $p\in
\partial f(E')$ is the image of a single point $x\in E':$
this is the reason to put the integrand function in (\ref{Mf}) equal
to 1, for every $p\in
\partial f(\Omega').$

Our purpose is to suggest a definition (see Theorem \ref{Monge
Ampere magari}) of the Monge--Amp\`ere measure of $u$ in the
Heisenberg context, on the analogy of the Euclidean framework,
using the H--normal map $\partial_H u$ of $u.$

\pn We know that if $u$ is a convex function in
$\Gamma^1(\Omega),$ then $\partial_H u(g)=\{\nabla_H u (g)\};$
however, it is unreasonable to require that $\nabla_H u:\Omega\to
V_1$ is a one--to--one map, since $\Omega\subset \H$ and $V_1$ is
essentially $\R^2.$ In other words, every point $v\in
\partial_H u(E)$ is the image of a set of points
$\Sigma_v^E\subset E.$ Therefore we need to replace the weight
\lq\lq 1\rq\rq\ in integral (\ref{Mf}) with a convenient weight.
For every $p\in V_1,$ the weight of $p$ will be the 2--dimensional
spherical Hausdorff measure of $\Sigma_v^E.$

In order to do this, we recall the following coarea formula proved
by Magnani in \cite{Ma2006}. We refer to \cite{Fe1969} for all the
relevant definitions about spherical Hausdorff measures.
\begin{theorem}
Let $F:\Omega\to\R^2$ be a Lipschitz map, where $\Omega\subset\H$
is a measurable set. Then, for every measurable function
$z:\Omega\to [0,\infty],$ the following formula holds
\begin{equation}\label{coarea}
 \int_\Omega z(g) J_HF(g)\, d g=\int_{\R^2}\left(
 \int_{F^{-1}(v)\cap \Omega}z(w)\, d
{\cal S}^2_\H(w)
  \right)\, dv,
\end{equation}
where $d{\cal S}^2_\H$ denotes the 2--dimensional spherical
Hausdorff measure.
 \end{theorem}

Let $u$ be a convex function in $\Gamma^2(\Omega),$ and consider
the function $F:\Omega\to\R^2$ defined by $F(g)=(Xu(g),Yu(g)).$
Clearly,
$$DF(g)=
    \left(
    \begin{array}{ccc}
    XXu(g)&YXu(g)&0\\
    XYu(g)&YYu(g)&0
    \end{array}
    \right).
$$
Standard computations give us that
$$
J_H F(g)=\texttt{\rm det} [\nabla_H^2 u](g)=\texttt{\rm det}
[\nabla_H^2 u(g)]^*+4(Tu(g))^2.
$$
 If we consider $z=1$ and $E\subset\Omega$ measurable,
by the formula (\ref{coarea}) we obtain\footnote{we identify
$\partial_H u(E)$ with a subset of $\R^2,$ as we did with $\nabla_H
u(g).$}
\begin{eqnarray*}
\int_E \left (\texttt{\rm det} [(\nabla_H^2 u)(g)]^* + 4
((Tu)(g))^2\right) d g &=& \int_E  J_HF(g)\, d g\\
 &=&\int_{\R^2}\left(
 \int_{F^{-1}(v)\cap E } d
{\cal S}^2_\H(w)
  \right)\, dv\\
 &=&
 \int_{\R^2} \left({\cal
S}^2_\H\bigl( (\nabla_H u)^{-1}(v) \cap E
\bigr) \right)\, d v\\
 &=&
 \int_{\partial_H u(E)}{\cal S}^2_\H\bigl(\Sigma^E_{v}
\bigr) \, d v,
\end{eqnarray*}
where, for every $v\in \R^2,$ the set  $\Sigma^E_{v}\subset \H$ is
defined by
 \begin{equation}\label{Sigma}
\Sigma^E_{v}= E\cap [(\nabla_H u)^{-1}(v)].
 \end{equation}

Taking into account the arguments above, we state the following
theorem where a possible definition for a Monge--Amp\`ere measure
associated to $u$ is provided.

\begin{theorem}\label{Monge Ampere magari}
Let $u\in{\Gamma}^{0,1}(\Omega),$ with $\Omega$ open in $\H.$ Let us
consider the function $\nu_u:{\cal E}\to [0,\infty]$ defined by
$$
\nu_u(E)=
 \int_{\partial_H u(E)}{\cal S}^2_\H\bigl(\Sigma^E_{v}
\bigr) \, d v,\qquad \forall E\in {\cal E},
$$
where $\Sigma^E_{v}$ is given in (\ref{Sigma}) and ${\cal
E}=\{E\subset\Omega:\ E\ \texttt{\rm and}\
\partial_H u(E)\ \texttt{\rm are Lebesgue measurable}\}.$    Then
\begin{enumerate}[i.]
\item $\nu_u$ is non negative and $\sigma$--additive;
\item if $u$
is convex and $u\in \Gamma^2(\Omega),$ then
$$
\nu_u(E)=\int_E \left (\texttt{\rm det} [(\nabla_H^2 u)(g)]^* + 4
((Tu)(g))^2\right) d g
$$ for every $E\in {\cal E}.$
\end{enumerate}
\end{theorem}
\psn We call $\nu,$ with an abuse of language, the Monge--Amp\`ere
measure associated to $u.$ Up to now, we are not able to prove
that ${\cal E}$ is a $\sigma$--algebra. Indeed, while it is quite
trivial that the numerable union of sets in ${\cal E}$ is still a
set in ${\cal E},$ it is not clear what happens about the
complement of a set in ${\cal E}.$ Notice that, for every
$E\in{\cal E},$ we have
$$
\partial_Hu(E^c)=\left(\partial_Hu(\Omega)\setminus
\partial_Hu(E)\right)\cup \left(\partial_Hu(\Omega\setminus E)\cap
\partial_Hu(E)\right).
$$
The main problem is to show that $\partial_Hu(\Omega\setminus
E)\cap
\partial_Hu(E)$ is
Lebesgue--measurable; notice that, in the Euclidean framework,
this set has null measure.
 \psn {\bf Proof:} First of all, let us
notice that for every $E\in{\cal E}$ and for every $v\in V_1,$ the
set $\Sigma^E_{v}$ is a Borel set; from the Borel regularity of
${\cal S}^2_\H,$ the set $\Sigma^{E}_{v}$ is ${\cal
S}^2_\H$--measurable.

Let us consider a sequence $\{E_i\}_{i=1}^\infty$ of disjoint sets
in ${\cal E}.$ It is straightforward that $ \Sigma^{\cup_i
E_i}_v=\cup_i \Sigma^{E_i}_v, $ and $\{\Sigma^{E_i}_v\}_i$ is a
family of disjoint subsets of $\Omega.$ We get
\begin{eqnarray*}
\nu_u(\cup_i E_i)&=&
 \int_{\R^2}{\cal S}^2_\H\bigl(\Sigma^{\cup_i E_i}_{v}
\bigr) \, d v\\
&=&
 \int_{\R^2}{\cal S}^2_\H\bigl(\cup_i \Sigma^{E_i}_{v}
\bigr) \, d v\\
&=&
 \int_{\R^2}\sum_i{\cal S}^2_\H\bigl( \Sigma^{E_i}_{v}
\bigr) \, d v\\
&=&
 \sum_i\int_{\R^2}{\cal S}^2_\H\bigl( \Sigma^{E_i}_{v} \bigr)\, d v\\
&=&
 \sum_i \nu_u( E_i).
\end{eqnarray*}
Hence, $\mu_u$ is $\sigma$--additive. Clearly, ii. is obvious for
previous computations  \QED

\vskip0.5truecm
\subsection{The Rockafellar function in $\H$}
\vskip0.5truecm

In the Euclidean framework, as well as in the more general Banach
setting, a convex function can be detected using its subdifferential
at every point via the Rockafellar function (for a new and recent
proof, see \cite{IvZl2008}). We are going to prove that a similar
integrability property is inherited by convex functions on the
Heisenberg group, where the H--subdifferential plays nearly the same
role. Indeed, the following result holds:

\begin{theorem}
Let $u:\H\to\R$ be a convex function. Then,
\begin{equation}\label{Rock function}
u(g)=u(g_0)+\sup_{\cal P} \left\{\sum_{i=0}^{n-1} \langle
p_i,\xi_1(g_{i+1})-\xi_1(g_i)\rangle\right\},
\end{equation}
were
$$ {\cal P}=\Bigl\{\{(g_i,p_i)\}_{i=0}^{n}: \ g_{i+1}\in H_{g_i},\ g_n=g,\
p_i\in\partial_Hu(g_i),\ n>0\Bigr\}.
$$
\end{theorem}
{\bf Proof:} Since if $g=g_0$ the result is trivial, we will assume
$g\neq g_0.$ From Theorem \ref{teorema mappa normale},
$\partial_Hu(g)\not =\emptyset,$ for every $g.$  By the definition
of H--subdifferential, for every sequence
$\{(g_i,p_i)\}_{i=0}^{n}\subset {\cal P}$ we have that
$$
u(g_{i+1})\ge u(g_i)+ \langle p_i,\xi_1(g_{i+1})-\xi_1(g_i)\rangle,
$$
for every $i,\ 0\le i\le n-1.$ Adding both sides of these
inequalities, we obtain
\begin{equation}\label{dim Rock 1}
u(g)\ge u(g_0)+ \sum_{i=0}^{n-1}  \langle
p_i,\xi_1(g_{i+1})-\xi_1(g_i)\rangle.
\end{equation}
Thus, the left--hand side of (\ref{Rock function}) is greater than
or equal to its right--hand side.

Since $g_{i+1}\in H_{g_i}$ if and only if $g_{i}\in H_{g_{i+1}},$
from $p_{i+1}\in
\partial _H u(g_{i+1}),$ we get that
$$
u(g_{i})\ge u(g_{i+1})+ \langle
p_{i+1},\xi_1(g_{i})-\xi_1(g_{i+1})\rangle.
$$
Then,
\begin{eqnarray*}
u(g_{i+1})-u(g_i)-\langle p_{i},\xi_1(g_{i+1})-\xi_1(g_{i})\rangle
&\le&
  \langle p_{i+1},\xi_1(g_{i+1})-\xi_1(g_{i})\rangle+\\
& &\qquad  -\langle
p_{i},\xi_1(g_{i+1})-\xi_1(g_{i})\rangle\\
&=&
 \langle p_{i+1}-p_i,\xi_1(g_{i+1})-\xi_1(g_{i})\rangle,
\end{eqnarray*}
for every $i,\ 0\le i\le n-1.$ Hence, taking into account (\ref{dim
Rock 1}), we obtain:
\begin{equation}\label{dim Rock 2}
0\le u(g)-u(g_0)-\sum_{i=0}^{n-1} \langle
p_{i},\xi_1(g_{i+1})-\xi_1(g_{i})\rangle \le
 \sum_{i=0}^{n-1} \langle p_{i+1}-p_i,\xi_1(g_{i+1})-\xi_1(g_{i})\rangle,
\end{equation}

In order to prove (\ref{Rock function}), we will show that, for
every $\epsilon>0,$ there exists a finite sequence in ${\cal P}$
such that
\begin{equation}\label{dim Rock 3}
 \sum_{i=0}^{n-1} \langle
 p_{i+1}-p_i,\xi_1(g_{i+1})-\xi_1(g_{i})\rangle
 \le \epsilon.
\end{equation}
Let us consider, first, the case $\xi_1(g_0)=\xi_1(g);$ this implies
that $g\not\in H_{g_0}.$ Choose $g'\in H_{g_0}$ and $g''\in
H_{g'}\cap H_g;$ from $g_0\not=g'\in H_{g_0},$ we have that
$\xi_1(g')\not=\xi_1(g)$ and hence $H_{g'}\cap H_{g}\not=\emptyset.$
Take $p_0\in \partial u_H(g_0),\ p'\in
\partial u_H(g'),\ p''\in \partial u_H(g'')$
and $p\in
\partial u_H(g).$ For every $\epsilon>0,$ denote by $N$ a positive integer such that
\begin{equation}\label{dim Rock 4}
\epsilon N\ge \langle p'-p_0,\xi_1(g')-\xi_1(g_{0})\rangle+
 \langle p''-p',\xi_1(g'')-\xi_1(g')\rangle+
 \langle  p-p'',\xi_1(g)-\xi_1(g'')\rangle
\end{equation}
We will single out a set of points $\{p_n\}_0^{3N}$ on the broken
line $[g_0,g']\cup[g',g'']\cup[g'',g].$ Indeed, let us consider a
particular sequence $\{(g_i,p_i)\}_{i=0}^{3N}$ in ${\cal P}$ defined
as follows:
\begin{itemize}
\item[{ i)\ }] for $i=1,\ldots N-1$, we pick out $g_i\in[g_0,g']\subset H_{g_0}$ such that $\xi_1(g_i)=(\xi_1(g')-\xi_1(g_0))i/N
+\xi_1(g_0);$

\item[{ii)\ }] we set $g_N=g'$ and $p_N=p';$

\item[{iii)\ }] for $i=N+1,\ldots 2N-1$, we pick out $g_i\in[g',g'']\subset\in H_{g'}$ such that $\xi_1 (g_i)=(\xi_1(g'')-\xi_1(g'))(i-N)/N
+\xi_1(g');$

\item[{iv)\ }] we set $g_{2N}=g''$ and $p_{2N}=p'';$

\item[{v)\ }] for $i=2N+1,\ldots 3N-1$, we pick out
$g_i\in[g'',g]\subset\in H_{g}$ such that $\xi_1
(g_i)=(\xi_1(g)-\xi_1(g''))(i-2N)/N +\xi_1(g'');$

\item[{vi)\ }] we set $g_{3N}=g$ and $p_{3N}=p:$

\item[{vii)\ }] for every $i$, with $1\le i\le 3N-1$ and $i$
different from $N$ and $2N,$ we choose $p_i\in
\partial u_H(g_i).$
\end{itemize}

Notice that $g_{i+1}\in H_{g_i},$ for every $i,\ 0 \le i\le 3N-1.$
From i)--vii) and (\ref{dim Rock 4}), we obtain
\begin{eqnarray*}
&&\sum_{i=0}^{3N-1} \langle
p_{i+1}-p_i,\xi_1(g_{i+1})-\xi_1(g_{i})\rangle=\\
&&\qquad= \left( \sum_{i=0}^{N-1}+ \sum_{i=N}^{2N-1}+
 \sum_{i=2N}^{3N-1}\right)
\langle
p_{i+1}-p_i,\xi_1(g_{i+1})-\xi_1(g_{i})\rangle\\
&&\qquad= \sum_{i=0}^{N-1} \langle
p_{i+1}-p_i,\xi_1(g')-\xi_1(g_{0})\rangle/N+ \sum_{i=N}^{2N-1}
\langle p_{i+1}-p_i,\xi_1(g'')-\xi_1(g')\rangle/N+\\
&&\qquad\qquad\qquad\qquad+
 \sum_{i=2N}^{3N-1}
\langle p_{i+1}-p_i,\xi_1(g)-\xi_1(g'')\rangle/N\\
&&\qquad= \frac{\langle
p_{N}-p_0,\xi_1(g')-\xi_1(g_{0})\rangle}{N}+\frac{
  \langle p_{2N}-p_N,\xi_1(g'')-\xi_1(g')\rangle}{N} +\\
&&\qquad\qquad\qquad\qquad\qquad\frac{\langle p_{3N}-p_{2N},\xi_1(g)-\xi_1(g'')\rangle}{N}\\
&&\qquad\le\epsilon.
\end{eqnarray*}
Hence (\ref{dim Rock 3}) holds.

If $g\notin H_{g_0}$ and $\xi_1(g)\not =\xi_1(g_0),$ we set $g'=g_0$
and choose $g''\in H_{g_0}\cap H_g.$ The proof is similar to the
previous case.

Finally, if $g\in H_{g_0},$ we set $g'=g_0$ and  $g''=g:$ again, the
proof follows the line of the previous case.
 \QED

\bibliography{cacapi}

\bibliographystyle{plain}
\end{document}